\documentclass[11pt]{article}
\long\def\kmdel#1{}
\long\def\kmcomment#1{}
\usepackage[a4paper,truedimen,scale={.80,.85},centering]{geometry}
\usepackage{latexsym,bm,setspace}
\usepackage{mathrsfs,amsmath,amssymb,amsthm,amscd,verbatim,ascmac}
\usepackage[titletoc,title]{appendix}
\usepackage[dvipdfm]{graphicx}
{
\theoremstyle{plain}
\newtheorem{thm}{Theorem}[section]

\newtheorem{kmProp}{Proposition}[section]

\newtheorem{kmRemark}{\textbf{Remark}}[section]
}

\newcommand{\rank}{{\rm rank}}
\newcommand{\Hd}{\mathcal{H}}
\newcommand{\frakC}{\mathfrak{C}}
\newcommand{\frakS}[1]{\mathfrak{S}_{#1}}
\newcommand{\frakg}{\mathfrak{g}}
\newcommand{\frakk}{\mathfrak{k}}
\newcommand{\fraksp}{\mathfrak{s}\mathfrak{p}}
\newcommand{\frakham}{\mathfrak{h}\mathfrak{a}\mathfrak{m}}
\newcommand{\Pkt}[2]{\{#1,#2\}}% 
\newcommand{\mR}{\ensuremath{\mathbb{R}}} %Use in Math-mode
\newcommand{\tz}{\tilde{z}}

\newcommand{\ds}{\displaystyle }
\newcommand{\myd}{d\,}
\newcommand{\mydz}{d\,}  % originally z means zero, but
\newcommand{\mydo}{d_{0}}  % originally o means one, but

\renewcommand{\tz}[2]{\ensuremath{\widehat{z}^{#1}_{#2}}}
\newcommand{\nz}[2]{{z}^{#1}_{#2}}
\newcommand{\kmqed}{\hfill \rule{1ex}{1.5ex}\par}

\newcommand{\CGF}[4]{\text{C}^{#1}_{\rm GF} ( {\frakham}_{#2}^{#3}, {\fraksp}(#2,\mR))_{#4}}
\newcommand{\HGF}[4]{\text{H}^{#1}_{\rm GF} ( {\frakham}_{#2}^{#3}, {\fraksp}(#2,\mR))_{#4}}
\newcommand{\cgf}[4]{\text{C}^{#1}_{\rm GF} ( {\frakham}_{#2}^{#3})_{#4}}
\newcommand{\hgf}[4]{\text{H}^{#1}_{\rm GF} ( {\frakham}_{#2}^{#3})_{#4}}
\newcommand{\km}{\ensuremath{ }} %Use in Math-mode
\newcommand{\pdel}{\partial } %Use in Math-mode
\allowdisplaybreaks

\begin{document}

\title{An affirmative answer to a conjecture for Metoki class\\
{\small To the memory of Professor Shoshichi Kobayashi at UC Berkeley}}

\author{Kentaro \textsc{Mikami}\thanks{ 
  Partially supported by Grant-in-Aid for 
  Scientific Research (C) of JSPS, No.26400063, No.23540067 
  and No.20540059}
}

%\classification{Primary 57R32, 57R17; Secondary 17B66}
% \keyword{Gel'fand-Kalinin-Fuks cohomology, Metoki class}

\label{startpage}

\maketitle

\begin{spacing}{1.25}

\abstract{In \cite{KOT:MORITA}, Kotschick and Morita showed 
 that the Gel'fand-Kalinin-Fuks class 
 in 
$\ds \HGF{7}{2}{ }{8}$  
 is decomposed as a product $\eta\wedge \omega $
  of some leaf cohomology class $\eta$ and 
  a transverse symplectic class $\omega$. 
  We show that the same formula holds 
  for Metoki class, which is a non-trivial 
  element in 
$\ds \HGF{9}{2}{ }{14}$.    
  The result was a conjecture 
  stated in \cite{KOT:MORITA}, where they studied 
  characteristic classes of symplectic 
  foliations due to Kontsevich. 
  Our proof depends on Gr\"obner Basis theory using  
  computer calculations.  
}

\section{Introduction}% 

This article is a superset of the author's preprint (\cite{KM:affirm}) with
the same title, which has been submitted to some Journal. As seen below,  in
order to give an affirmative proof for D.~Kotschick and S.~Morita's
conjecture, we have had to handle the enormous data of the concrete bases of
the cochain complexes, matrix representations of coboundary operators, and
proceed with many symbol calculus calculations by computer. These data are
available at the author's private URL
\text{www.math.akita-u.ac.jp/\~{}mikami/Conj4MetokiClass/}.  The Journal
requires that those data should be  stored in a ``perpetual'' place.  Thus,
the author appended these data to the original paper as appendices, then the
page number of the enlarged paper reached 1167.   By much help of arxiv.org,
the author stored  the data files of plain text format to the subdirectory
\texttt{anc/} as ancillary files, and is trying to upload the package to
arxiv.org.   

Let $\ds {\frak X}(M)$ be the Lie algebra of smooth vector fields of a
smooth manifold $M$.  $\ds\text{H}_{c}^{\bullet}(  {\frak X}(M)) $ is the
Lie algebra cohomology groups, where the subscript $c$ means that each
cochain is required to be {\em continuous}. The cohomology group
$\ds\text{H}_{c}^{\bullet}(  {\frak X}(M)) $ is often written as
$\ds\text{H}_{\text{GF}}^{\bullet}(M)$ and is called the Gel'fand-Fuks
cohomology group of $M$.    
It is known that if $M$ is of finite-type (i.e., $M$ has a finite cover),
then 
$\text{H}_{\text{GF}}^{\bullet}(M)$ is finite dimensional. 

%%%%%%%%%%%%%%%%%%%%%%%%%%%%%%%%%%%%%%%%%%%%%%%%%%%%%%%%%%%%%%%%%%% 
Let 
$\ds {\frak a}_{n}$ denote the Lie algebra of 
formal vector fields on $\ds \mR^{n}$ which is given by 
$\ds 
\mR[[x_1,\ldots, x_n]]\langle 
\pdel/\pdel x_1,\ldots,  \pdel/\pdel x_n\rangle $. 
Thus an element of 
$\ds {\frak a}_{n}$ is a vector field with coefficients of formal
power series in coordinate functions. Then  
$\ds\text{H}_{c}^{\bullet}(  {\frak a}_{n} ) \cong 
\text{H}_{\text{GF}}^{\bullet}( \mR^ {n}) $ holds and so 
$\ds\dim\text{H}_{c}^{\bullet}(  {\frak a}_{n} ) <\infty$.  

%%%%%%%%%%%%%%%%%%%%%%%%%%%%%%%%%%%%%%%%%%%%%%%%%%%%%%%%%%%%%%%%%%%

Let 
$\ds{\frak v}_{n}$ be a subalgebra of 
${\frak a}_{n}$, consisting of 
the volume preserving formal vector fields on $\ds\mR^{n}$, 
and
$\ds
\frakham_{2n}$  a subalgebra of 
${\frak a}_{2n}$, consisting of   
formal Hamiltonian vector fields on $\ds\mR^{2n}$. 
Then, the next question is still open: 
Is $\ds \dim \text{H}_{c}^{\bullet} ({\frak v}_{n})$  infinite or  
is 
$\ds\dim \text{H}_{c}^{\bullet} (\frakham_{2n})$ infinite?  
%%%%%%%%%%%%%%%%%%%%%%%%%%%%%%%%%%%%%

\medskip

There is a notion of {\em weight} for cochains of $\ds
\frakham_{2n}$. Since the weight is preserved by the coboundary
operator, there is the cohomology subgroup corresponding to each
weight (cf.\ \S \ref{wgt}(I-3)).  
In \cite{MR0312531}, %   Gel'fand-Kalinin-Fuks got 
(1) For the weight $w \le 0$, the structure of 
$\ds \text{H}_{c}^{\bullet} (\frakham_{2n}, \fraksp(2n,\mR))_{w}$ is
completely determined, and   

(2) When $n=1$ and $w>0$, the next holds true:
\begin{align*}
\text{H}_{c}^{\bullet} (\frakham_{2}, \fraksp(2,\mR))_{w} &= 0 \quad (w=1,2,\ldots,7) \\   
\text{H}_{c}^{m} (\frakham_{2}, \fraksp(2,\mR))_{8} &= 
        \begin{cases}\ \mR & m=7\\ \ 
0 & \text{otherwise} \end{cases}
\end{align*} 
The generator of 
$\ds \text{H}_{c}^{7} (\frakham_{2}, \fraksp(2,\mR))_{8}$ is called   
the Gel'fand-Kalinin-Fuks class.   
Hereafter, we use the notation 
$\ds\HGF{\bullet}{2}{ }{w} $ instead of 
$\ds\text{H}_{c}^{\bullet} (\frakham_{2}, \fraksp(2,\mR))_{w}$. 

There is a homomorphism from $\ds \HGF{\bullet}{2}{ }{ }$ into
$\ds\text{H}^{\bullet}( B\Gamma_{2}^{symp})$, where $\ds
\Gamma_{2}^{symp}$ is the groupoid of germs of local
diffeomorphisms of $\ds\mR^2$ preserving the symplectic structure of
$\ds\mR^2$.  
It is still not known if the image of the Gel'fand-Kalinin-Fuks class
by the homomorphism is trivial in $\ds\text{H}^{7} (
B\Gamma_{2}^{symp})$ or not (cf.\cite{B:H:MR0307250}).    

The next non-trivial result in succession to the
Gel'fand-Kalinin-Fuks class is   $\ds  \HGF{9}{2}{ }{14} \cong
\mR$, which was shown by S.~Metoki (\cite{metoki:shinya}) in 1999.
He is interested in the volume preserving formal vector fields,
when $n=1$ both $\ds \frakham_{2}$ and $ {\frak v}_{2}$ are the
same.  

Let $\ds {\cal F}$ be a foliation on a manifold $M$. We have the
foliated cohomology defined by $\ds \text{H}_{\cal F}^{\bullet}
(M,\mR) := \text{H}^{\bullet}( \Omega_{\cal F})$ where $\ds
\Omega_{\cal F} = \Omega(M)/\text{I}({\cal F})$, $\ds\Omega(M)$ is
the exterior algebra of differential forms on $M$, and
$\ds\text{I}({\cal F})$ is the ideal generated by $\ds\{ \sigma
\in \Omega^{1}(M) \mid \langle \sigma, \text{T}{\cal F}\rangle =0
\}$.   

M.~Kontsevich (\cite{Kont:RW}) showed that if $\ds {\cal F}$ is a
codimension $2n$ foliation endowed with a symplectic form 
$\omega$ in the transverse direction, then there is a commutative
diagram: $$\begin{CD} @. \text{H}_{\cal F}(M,\mR) @>\omega^{n}\wedge>>
        \text{H}_{\text{ DR}}^{\bullet+2n}(M,\mR)  @.  \\ @.
        @AAA   @AAA \\ @.   \HGF{\bullet}{2n}{0}{ } @>\omega^{n}\wedge>>
        \HGF{\bullet+2n}{2n}{ }{ } @.  \end{CD} $$ where   $\ds
        \frakham_{2n}^{0}$ is the Lie subalgebra of the
        Hamiltonian vector fields of the formal polynomial vanishing 
        at the origin of $\ds  \mR^{2n}$.  
\kmcomment{ About the recent developments like as (commutative)
graph (co)homology and so on,  concerning to this new view point
by M.~Kontsevich, we refer to Professor S.~Morita's special
lecture series of Encounter with Math at Chuo University.  }

D.~Kotschick and S.~Morita (\cite{KOT:MORITA}) 
determined the space
$\ds  \HGF{\bullet}{2}{0}{w} $  
for $w\le 10$, and    
concerning Kontsevich homomorphism in the case of $n=1$, 
\kmcomment{
		$$\ds \HGF{\bullet}{2n}{}{w}
\mathop{\longrightarrow}^{\wedge \omega^{n}}
\HGF{\bullet+2n}{\bullet}{2n}{1}{w-2n}$$ 
} 
they showed the following, as well as the non-triviality of
Kontsevich homomorphism:   

\begin{thm}[\cite{KOT:MORITA}]\label{thm:morita}\rm 
There is a unique element $\ds  \eta \in
\HGF{5}{2}{0}{10} \cong \mR$ such that 
$$ \text{Gel'fand-Kalinin-Fuks class } = \eta \wedge \omega 
\in \HGF{7}{2}{ }{8}  $$
where $\omega$ is the cochain associated with the 
linear symplectic form of $\ds  \mR^2$. 
\end{thm}

Further they stated that it is highly likely that the same thing is true also for 
$\ds \text{Metoki class}\in  \HGF{9}{2}{ }{14}$.  That is, 
there should 
exist an element 
$\ds  \eta' \in \HGF{7}{2}{0}{16} 
$ such that 
$$ \text{Metoki class} = \eta' \wedge \omega 
\in \HGF{9}{2}{ }{14}\ . $$ 
In the same line of 
D.~Kotschick and S.~Morita (\cite{KOT:MORITA}), 
we determined  $\ds\HGF{\bullet}{2}{0}{w} $ 
for $w\le 20$ in \cite{M:N:K}.  
In this paper, making use of information in \cite{M:N:K}, 
we will show the following theorem.

\medskip

\begin{thm}\label{thm:mikami}\rm 
$\ds \HGF{9}{2}{ }{14}$ and 
$\ds\HGF{7}{2}{0}{16}$ are both 1-dimensional and 
the map of wedging symplectic cocycle 
$$ \omega\wedge : \HGF{7}{2}{0}{16} \longrightarrow \HGF{9}{2}{ }{14} $$
is an isomorphism.   
$$ \omega\wedge : \HGF{7}{2}{0}{16} \longrightarrow \HGF{9}{2}{ }{14} $$
is an isomorphism.   
Thus, 
there is a unique element $\ds  \eta' \in
%\text{H}_{\text{GF}}^{7}(\frakham_{2}^{0})_{16}^{Sp}
\HGF{7}{2}{0}{16} 
\cong \mR$ such that 
$$ \text{Metoki class } = \eta' \wedge \omega 
\in \HGF{9}{2}{ }{14} $$
where Metoki class is the generator of 
$\ds  \HGF{9}{2}{ }{14}$.  
\kmcomment{
Shinya Metoki's Doctoral Thesis
\cite{metoki:shinya}, 
which we could get from Prof.~T.~Tsuboi's Web Page at the University of
Tokyo.    
}
\end{thm} 
%\vspace{3cm}

%\input body-GB.tex

\section{Preliminaries}
Generalities concerning 
the (relative) Gel'fand-Fuks
cohomologies and symplectic formalism are found in  
Mikami-Nakae-Kodama's preprint (\cite{M:N:K}).   Here we review the concept
of cochain complex of our Lie algebras and the symplectic action on relative
complex and also the description of coboundary operators for further
calculations. 
Although the space we are concerned with in this paper is $\ds 
\mR^2$,  
we review the notions 
on the general linear symplectic space 
$\ds  \mR^{2n}$, 
and fix notations we use hereafter.

% \subsection{$2n$-dimensional symplectic space $\ds  \mR^{2n}$} 
\subsection{Symplectic space $\ds\mR^{2n}$} \label{wgt}
We fix a linear
symplectic manifold $\ds  (\mR^{2n}, \omega)$ with the standard
variables $x_1,x_2,\ldots,x_{2n}$. 
Let $\ds \Hd_f$ denote the Hamiltonian vector field of $f$.    
Recalling the formula $\ds  [\Hd_f, \Hd_g] = -
\Hd_{\Pkt{f}{g}}$ for Hamiltonian vector fields, 
we identify each formal Hamiltonian vector
field with its potential polynomial function up to the constant term and the Lie
bracket of vector fields 
with the Poisson bracket on polynomial functions, where 
We denote by $\ds  \frakS{p}$     
the dual space of homogeneous polynomials of degree $p$.  
Then     
\begin{align*}
        \ds  \frakham_{2n}^{ } =& \left(\mathop{\oplus}_{p=1}^{\infty}
\frakS{p}^{*}\right)^{\wedge} \quad \text{is a Lie algebra }\\
\noalign{and}
\frakham_{2n}^{0} =& \left(\mathop{\oplus}_{p=2}^{\infty}
\frakS{p}^{*}\right)^{\wedge} \quad \text{is a subalgebra of } \frakham_{2n}^{ },  
\end{align*} 
where $\ds \left(\phantom{MM}\right)^{\wedge}$ means the completion 
using the Krull topology. 

Using the above notation we have the following:

(I-1) $m$-th cochain complexes of $\ds\frakham_{2n}$ and
$\ds\frakham_{2n}^{0}$  are given by  
\begin{align*}
C_{GF}^{m}( \frakham_{2n}^{ }) =& \mathop{\oplus}_{k_1+k_2+\cdots =m}
\Lambda^{k_1}\frakS{1} \otimes
\Lambda^{k_2}\frakS{2} \otimes
\Lambda^{k_3}\frakS{3} \otimes \cdots \\
\noalign{and $\ds\frakS{1}$ is the dual space of constant vector fields}
C_{GF}^{m}( \frakham_{2n}^{0}) =& \mathop{\oplus}_{k_2+k_3+ \cdots =m}
\Lambda^{k_2}\frakS{2} \otimes
\Lambda^{k_3}\frakS{3} \otimes
\Lambda^{k_4}\frakS{4} \otimes \cdots \ . 
\end{align*}
%The difference above is only $k_1=0$.  

(I-2) 
The coboundary operator $\mydz$ on
$ \ds  \cgf{\bullet}{2n}{}{}$ 
is 
defined by 
$$ (\mydz  \sigma) (f_0,f_1,\ldots,f_m) = 
\sum_{k<\ell } (-1)^{k+\ell } \sigma(
\Pkt{f_k}{f_\ell }, \ldots, \widehat{f_k},\ldots , \widehat{f_\ell },
\ldots) \qquad f_i \in \frakham_{2n}$$ 
for each $m$-cochain 
$ \sigma\in \ds  \cgf{m}{2n}{}{}$. 

And the coboundary operator $\mydo$ on
$ \ds  \cgf{\bullet}{2n}{0}{}$ 
is 
defined by 
$$ (\mydo  \sigma) (f_0,f_1,\ldots,f_m) = 
\sum_{k<\ell } (-1)^{k+\ell } \sigma(
\Pkt{f_k}{f_\ell }, \ldots, \widehat{f_k},\ldots , \widehat{f_\ell },
\ldots) \qquad f_i \in \frakham_{2n}^{0}$$ 
for each $m$-cochain 
$ \sigma\in \ds  \cgf{m}{2n}{0}{}$.

We will study the difference between two coboundary operators $\mydz$
and $\mydo$ in subsection \S\ref{subsec::coboundary}.  

(I-3) There is a notion of \textbf{weight} for cochains
(cf.\cite{KOT:MORITA}).   
For each non-zero cochain 
$$\sigma \in 
\Lambda^{k_1}\frakS{1} \otimes
\Lambda^{k_2}\frakS{2} \otimes
\Lambda^{k_3}\frakS{3} \otimes \cdots \otimes 
\Lambda^{k_\ell}\frakS{\ell} 
$$ 
its weight is given by $$(1-2)k_1 + (2-2)k_2 + (3-2)k_3 +\cdots +
(\ell-2)k_{\ell} = \sum_{i=1}^{\ell} (i-2) k_i\ . $$

The weight of a cochain is preserved by the coboundary operator, and we
can decompose  each cochain
complex by way of weights and get Gel'fand-Fuks cohomologies with a discrete
parameter, namely with \textit{weight} $w$  like as 
$$ \cgf{m}{2n}{j-1}{w}\qquad\text{and}\qquad  
 \hgf{m}{2n}{j-1}{w} \quad  
 \text{for } j=0,1 $$  
where 
$\ds \frakham_{2n}^{-1}$ means the space $\ds \frakham_{2n}$.

In both cases, for given degree $m$ and weight $w$, 
we consider the sequences $(k_1,k_2,k_3,\ldots)$
of nonnegative integers 
with
\begin{equation} \sum_{j=1}^{\infty} k_j = m \qquad\text{and}\qquad  
 \sum_{j=1}^{\infty}(j-2) k_j = w\ . 
 \label{degree:org:weight}
\end{equation}
Readers may be anxious about the contribution of $\ds  k_2$ or $
\ds  
k_1$. In fact,  
there is a dimensional restriction for each $\ds  k_j$ with  
$\ds  0 \le k_j \le \dim \frakS{j} = \dfrac{ (j+2n-1)!}{j!
(2n-1)!}$.

From those two relations in (\ref{degree:org:weight}), we have
\begin{equation} \label{degree:weight}
 \sum_{j=1}^{\infty} k_j =  m  
 \quad\text{and}\quad 
 \sum_{j=1}^{\infty}j k_j = w +2 m  \ . 
 \end{equation} 
This means our sequences correspond to all partitions of $w+2m$ of
length $m$,  or in other words, 
to the Young diagrams with $w+2m$ cells of length $m$ (cf.\cite{M:N:K}).  
Furthermore, we require dimensional restrictions, and
$\ds  k_1=0$ when $\ds  \frakham_{2n}^{0}$.

\subsection{Symplectic action and the relative cohomologies} 
%%\begin{frame}[allowframebeaks]{Relative cohomologies} 
We denote the natural action of 
the Lie group $K = Sp(2n,\mR)$ on $\ds  \mR^{2n}$ by $\varphi_a$ for $a\in K$, 
i.e., 
$\varphi_a( \mathbf{x}) = a \mathbf{x} $ as the multiplication of
matrices.  The action leaves $\omega$ invariant by definition, and 
we see that 
$\ds  (\varphi_a)_{*} (\Hd_f) = \Hd_{ f\circ \varphi_{a^{-1}}}$ 
for
each function $f$ on $\ds  \mR^{2n}$ and $a\in K$.
Let $\ds  \frakk  = \fraksp(2n,\mR)$ be the Lie algebra of $K$.     
We denote 
the fundamental vector field on $\ds  \mR^{2n}$
of $K$ by $\ds  \xi_{\mR^{2n}}$ for $\xi\in\frakk$. 
The equivariant (co-)momentum mapping of symplectic action of $K$ is given by  

$$ \hat{J}(\xi) \mathbf{x}  = -\frac{1}{2} {}^t \mathbf{x} 
\begin{bmatrix}
        \Pkt{x_1}{x_1} & \ldots & \Pkt{x_1}{x_{2n}} \\
        \vdots  & \ldots & \vdots  \\
        \Pkt{x_{2n}}{x_1} & \ldots & \Pkt{x_{2n}}{x_{2n}} \\
\end{bmatrix} %\left( \Pkt{x_i}{x_j} \right) 
\xi \mathbf{x} $$ 
where $\mathbf{x}$ is the natural coordinate of $\ds  \mR^{2n}$
as column vector,  
$\ds  {}^t \mathbf{x}$ means the transposed row vector of $\ds\mathbf{x}$, 
$\ds \Pkt{x_i}{x_j}$ is the Poisson bracket of $i$-th
and $j$-th components of $\mathbf{x}$ with respect to $\omega$, 
and  
$\xi\in \frakk$.   $\ds  \hat{J}$ is a Lie algebra monomorphism
from the Lie algebra $\fraksp(2n,\mR)$ into the Lie algebra
$\ds  
C^{\infty}(\mR^{2n})$ with the Poisson bracket. We stress that 
$\ds  
\hat{J}(\xi) $ is a degree 2 homogeneous polynomial function on $\ds 
\mR^{2n}$ for $\xi \ne 0$.  
%
%\begin{pmatrix} O & I \\-I & O \end{pmatrix} 
%
%
The Hamilton potential of the bracket $\ds 
[\xi_{\mR^{2n}}, \Hd _f]$ is given by $\ds  - \Pkt{\hat{J}(\xi)}{f}$,
because of $\ds  [\xi_{\mR^{2n}}, \Hd_f] = [\Hd_{\hat{J}(\xi)}, \Hd_f] 
= - \Hd_{ \Pkt{ \hat{J}(\xi)}{f} } 
$.  This means that 
$\frakk$ is regarded as a subalgebra of 
$\mathfrak g = \frakham_{2n}^{ } $
or 
$\mathfrak g = \frakham_{2n}^{0} $ through the equivariant momentum
mapping $J$.    

Define the relative cochain group 
$\ds  \text{C}^{m}( {\mathfrak g}, {\mathfrak k} ) $ by 
$$\ds  \text{C}^{m}( {\mathfrak g}, {\mathfrak k} ) 
= \{ \sigma \in
\text{C}^{m}({\mathfrak g}) \mid i_X\sigma =0, i_X \myd \sigma =0\quad
(\forall  X\in {\mathfrak k})\}
\quad (m=0,1,\ldots)\ .
$$
Then  $\ds\left( \text{C}^{\bullet}(\frakg, \frakk), d\right)$ becomes a cochain
complex,   
and we get the relative cohomology groups $\ds 
\text{H}^{m}({\frakg}, {\frakk})$.  
Let $K$ be a Lie group of $\mathfrak k$. Then we also define 
$$\ds  \text{C}^{m}( {\mathfrak g},  K ) = \{ \sigma \in
\text{C}^{m}({\mathfrak g}) \mid i_X\sigma =0\  
(\forall  X\in {\mathfrak k}), Ad_k^{*} \sigma = \sigma \ (\forall k\in K)
\}$$ and we get the relative cohomology groups $\ds 
\text{H}^{m}({\mathfrak g}, K)$. If $K$ is connected, these two cochain
groups are identical.  
If $K$ is a closed subgroup of $G$, then it can be seen $\text{C}^{\bullet}({\mathfrak g}, K)
= \Lambda ^{\bullet} ( G/K) ^{G} $ (the exterior algebra of
$G$-invariant differential forms on $G/K$).

Since the space of degree 2 homogeneous polynomials, $\ds  \frakS{2}^{*}$ is
spanned by the image of momentum mapping $\ds  \hat{J}$ of
$Sp(2n,\mR)$, we see that 
\begin{kmProp}[\cite{M:N:K}]\label{m:n:k:decomp} \rm 
For each cochain $\sigma$, $i_{\xi} \sigma=0$ 
($\forall \xi\in\fraksp(2n,\mR)$) implies $k_2=0$, 
and the other condition $i_{\xi} d\sigma=0$ is equivalent to 
$L_{\xi} \sigma=0$ ($\forall \xi\in\fraksp(2n,\mR)$).  
Thus we see  for $j=0,1$ 
\begin{align*}
&
 \CGF{\bullet}{2n}{j-1 }{w} =  
 %\sum_{(k_1,k_2,k_3,\ldots) \text{ with Cond}_{j}} 
 \sum_{\text{Cond}_{j}} 
        \left( \Lambda^{k_1} \frakS{1} \otimes 
 \Lambda^{k_2} \frakS{2} \otimes 
 \Lambda^{k_3} \frakS{3} \otimes \cdots\right)^{triv}  
\end{align*} 
where $\ds \left(\phantom{MM}\right)^{triv}$ means the direct sum of the 
(underlying) subspaces of the trivial representations.   
$\ds \text{Cond}_0$ consists of the conditions 
(\ref{degree:weight}) in the preceding subsection,  
$\ds k_2=0$, and the dimensional
restrictions.     $\ds \text{Cond}_1$ consists of 
$\ds \text{Cond}_0$ and $\ds k_1=0$.    
\end{kmProp}

As already explained in 
\cite{KOT:MORITA}, if the weight $w$ is odd 
$\ds \HGF{\bullet}{2n}{j-1}{w} = 0$ for $j=0,1$. Thus, we
have only to deal with even weights.  

\begin{kmRemark}\rm 
There is a notion of \textbf{type} $N$ for cochains in 
\cite{metoki:shinya}.  The weight $w$ and type $N$ are related by
$w = 2 N $.   

There is a general method to decompose $\ds\Lambda^{p}\frakS{q}$ into the
irreducible subspaces for a given $Sp(2n,\mR)$-representation,  by  
getting the maximal vectors which are invariant by the maximal unipotent subgroup of $Sp(2n,\mR)$.

Concerning the decomposition of the tensor product, we have the
Clebsch-Gordan rule when $n=1$.  (For $n=2$, Littlewood-Richardson
rule is used in \cite{Mik:Nak}, and the crystal base theory is
used in \cite{KM:D6} when $n=3$.)  
\end{kmRemark}

%%\end{frame}
%\begin{frame}[allowframebreaks]{Metoki's way}

%\newcommand{\mytime}{\text{time}}
\newcommand{\mytime}{\text{}}

\subsection{Coboundary operators} \label{subsec::coboundary}
We want to distinguish the coboundary operators where they act.  
By $\mydz$, we mean the coboundary operator which acts on $\ds  \CGF{\bullet}{2n}{ }{w}$ and 
by $\mydo$, the one on $\ds  \CGF{\bullet}{2n}{0}{w}$.   

Let $\omega$ be the 2-cochain defined by the linear symplectic form of
$\ds  \mR^{2n}$, we see that $$\ds  \omega \in 
\CGF{2}{2n}{ }{(-2)} {\setminus}  
\CGF{2}{2n}{0}{(-2)}$$ and   
$\ds  \omega^{n} \in 
\CGF{2n}{2n}{ }{(-2n)}$.   

\begin{kmProp}\label{cor:one}
The linear map $$\ds  \omega^n \wedge : 
\CGF{\bullet}{2n}{0}{w} \longrightarrow 
\CGF{\bullet+2n}{2n}{ }{w-2n}$$ satisfies 
$$ \mydz ( \omega^n \wedge \sigma ) = \omega^n \wedge \mydo (\sigma)$$ and 
so the next diagram is commutative 
\begin{equation}\label{comm:diagram}
\begin{CD}
@. \CGF{\bullet-1+2n}{2n}{ }{w-2n}@>\mydz>>\CGF{\bullet+2n}{2n}{ }{w-2n} @.  \\
@.   @A\omega^{n}\wedge AA   @AA\omega^{n}\wedge A \\
 @.   \CGF{\bullet-1}{2n}{0}{w} @>\mydo>> \CGF{\bullet}{2n}{0}{w} @.
\end{CD}
\end{equation} 
Thus we have a linear 
 map $$\ds  \omega^n \wedge : 
\HGF{\bullet}{2n}{0}{w} \longrightarrow 
\HGF{\bullet+2n}{2n}{ }{w-2n}$$ 
naturally.  
This induced map is trivial if and only if
\begin{equation} %\label{triv:cond}
\omega^{n}\wedge \ker(\mydo \text{ on }
\CGF{\bullet}{2n}{0}{w} 
%\rightarrow \CGF{\bullet+1}{2n}{0}{w} 
) \subset \mydz\left( 
\CGF{\bullet-1+2n}{2n}{ }{w-2n}\right)
\label{simple:key:issue}
\end{equation}
\end{kmProp} 
Proof: 
We have 
$\ds  \mydz( \omega ) = 0$. This is a requirement for a symplectic form.  
For each $\ds  \sigma \in \cgf{\bullet}{2}{0}{}$, 
we already know that $\ds  \omega^n \wedge (\mydz(\sigma) -
\mydo(\sigma)) = 0$ and now we see that  
$\ds  \omega^n \wedge \mydz(\sigma) = \mydz ( \omega^n \wedge
\sigma) $ because of $\ds  \mydz \omega=0$. 
The above states that only the diagram (\ref{comm:diagram}) is commutative. 
Then we have 

(i) If $\ds  \mydo(\sigma)=0$, then 
$\ds  \mydz(\omega^n \wedge \sigma)=0$,   
namely, 
$\ds  \omega^n \wedge \ker(\mydo) \subset \ker(\mydz)$.   

(ii) If  
$\ds  \sigma,\tau \in \CGF{\bullet}{2n}{0}{}$ satisfy 
$\ds  \mydo(\sigma) = 0 = \mydo(\tau)$ and 
$\ds  \sigma -\tau = \mydo(\rho)$,  then 
$\ds \omega^n\wedge\sigma-\omega^n\wedge\tau=\omega^n
\wedge \mydo(\rho) = \mydz ( \omega^n \wedge \rho)$.  This means that 
the wedge product by 
$\ds  \omega^n$ induces a well-defined linear map 
$$\ds  \omega^n : \HGF{\bullet}{2n}{0}{} \longrightarrow 
\HGF{\bullet+2n}{2n}{ }{}\quad \text{by}\quad 
\sigma \mapsto \omega^n \wedge \sigma \ . 
$$   

(iii) From (i), we see that the map is trivial if and only if 
$\ds  \omega^n \wedge \ker(\mydo) \subset \mydz(
\CGF{\bullet+2n-1}{2n}{ }{w-2n} )$.   

\kmqed

%\section{Proof for the conjecture}
%\section{The symplectic 2-plane $\ds  \mR^2$}
\section{Symplectic 2-plane}\label{sec:three}
In this section, we deal with the symplectic 2-plane 
$\ds\mR^2$. 
We study the difference between two coboundary operators $\mydz$ and
$\mydo$.  Since $\ds  \dim \frakS{1}=2$, the
domain of definition 
of $\mydz$ is given by    
$$ \ds  \cgf{\bullet}{2}{ }{} = 
 \ds  \cgf{\bullet}{2}{0}{}
 \oplus \left(\frakS{1} \otimes 
 \ds  \cgf{\bullet-1}{2}{0}{}
 \right)
 \oplus \left(\Lambda^2 \frakS{1} \otimes 
 \ds  \cgf{\bullet-2}{2}{0}{} 
 \right)\, . 
 $$ 
Let $x$, $y$ be a global Darboux coordinate satisfying $\Pkt{x}{y}=1$.  
We denote 
by 
$\ds  \nz{r}{R} $,
the dual element of the 
polynomial $\ds  \tz{r}{R} = 
\dfrac{x^r}{r!} \dfrac{ y^{R-r}}{(R-r)!}$ where 
$R>0$ and $0\le r \le R$.

The two coboundary operators $\mydz$, $\mydo$ in those bases, are 
\begin{align*}
\mydz \nz{r}{R} &= - \frac{1}{2} \sum_{A+B=2+R} \langle \nz{r}{R},
\Pkt{\tz{a}{A}}{\tz{b}{B}}\rangle  
\nz{a}{A} \wedge \nz{b}{B}
\end{align*}
where $A>0$, $B>0$, $ 0\le a\le A$,    $0\le b\le B$.   
\begin{align*}
\mydo \nz{r}{R} &= - \frac{1}{2} \sum_{A+B=2+R}\langle \nz{r}{R},
\Pkt{\tz{a}{A}}{\tz{b}{B}}\rangle  
\nz{a}{A} \wedge \nz{b}{B}
\end{align*}
where $A>1$, $B>1$, $0\le a\le A$,    $0\le b\le B$.   
Thus, the difference between $\ds  \mydz $
and $\ds  \mydo $ for a 1-cochain can be written as follows.   
\begin{align*} 
\mydz \nz{r}{R} 
=&  \mydo \nz{r}{R} - \sum_{0\le a\le 1, 0\le b\le 1+R} \langle \nz{r}{R}, 
\Pkt{ \tz{a}{1}}{ \tz{b}{1+R}} \rangle 
 \nz{a}{1} \wedge  \nz{b}{1+R}
 \\
=&  \mydo \nz{r}{R} 
+ \nz{0}{1}  \wedge  \nz{1+r}{1+R}
- \nz{1}{1}  \wedge  \nz{r}{1+R} 
= 
 \mydo \nz{r}{R} 
+ \begin{vmatrix}
\nz{0}{1}  & \nz{r}{1+R} \\
\nz{1}{1}  & \nz{1+r}{1+R} 
\end{vmatrix}  
\end{align*}
where
$\ds  
\begin{vmatrix}
\nz{0}{1}  & \nz{r}{1+R} \\
\nz{1}{1}  & \nz{1+r}{1+R} 
\end{vmatrix}  
$ is the determinant of the (2,2)-matrix whose multiplication is the
wedge product.

We may assume that $\ds  \mydo \nz{r}{1} = 0$.  
The 2-cochain 
$\omega$ which comes from the symplectic
structure, is written as $\ds  
\omega = \nz{0}{1} \wedge \nz{1}{1}$ in our notation and we see directly that 
\begin{align*} \mydz \omega =& \mydz (\nz{0}{1} \wedge \nz{1}{1}) 
=  
\begin{vmatrix}
\nz{0}{1} & \nz{0}{2} \\
\nz{1}{1} & \nz{1}{2} \end{vmatrix}
\wedge \nz{1}{1} - 
\nz{0}{1} \wedge  
\begin{vmatrix}
\nz{0}{1} & \nz{1}{2} \\
\nz{1}{1} & \nz{2}{2} \end{vmatrix} 
=  
- \begin{vmatrix}
\omega  & \nz{0}{2} \\
0       & \nz{1}{2} \end{vmatrix}
 - 
\begin{vmatrix}
0       & \nz{1}{2} \\
\omega  & \nz{2}{2} \end{vmatrix} 
= 0 \, . 
\end{align*}
But, $\omega$ is not $\mydz$-exact because  
$\ds  \Pkt{\tz{a}{1}}{\tz{b}{1}}=\text{constant}$.

\section{Proof of Theorem}  
In this section, we give a proof for  
Theorem \ref{thm:mikami} which asserts that 
$$ \omega\wedge :
\HGF{7}{2}{0}{16} \rightarrow  
\HGF{9}{2}{ }{14} $$ 
is an isomorphism.    Since we know that the source and the target spaces
are both 1-dimensional, it is enough to show the map $\omega\wedge$ is
non-trivial.  For that purpose, we make use of (\ref{simple:key:issue}) of
Proposition \ref{cor:one}.

We have information about 
$\ds  \CGF{\bullet}{2}{0}{w}$ ($w=12,14,16,18,20$) (cf. \cite{M:N:K}).  
We show the result of weight =16 in the table below. 
In the table 
$\ds  C^k$ is     
$\ds  \CGF{k}{2}{0}{16}$,  
$\dim$ is $\dim C^k$, and   
$\rank$ is the rank of $ \ds  \mydo : C^k \longrightarrow
C^{1+k}$.  

\begin{center}
\tabcolsep= 4pt
\begin{tabular}{|c|*{14}{c}c|}
\hline
$\ds  \frakham_{2}^{0}$, w=16
 & $\ds  
 \mathbf{0}$ & $\rightarrow$ &
$\ds  C^{3} $& $ \rightarrow $ &
$\ds  C^{4} $& $ \rightarrow $ &
$\ds  C^{5} $& $ \rightarrow $ &
$\ds  C^{6} $& $ \rightarrow $ &
$\ds  C^{7} $& $ \rightarrow $ &
$\ds  C^{8} $& $ \rightarrow $ &
$\mathbf{0}$ \\
\hline
$\dim$  &   && 12 && 61 && 126 && 147 && 95&& 24 && \\
rank & & 0 && 12 && 49 && 77 && 70 && 24 && 0 &\\
%$\dim$(ker) &   && 0 && 12 && 49 && 77 && 71&& 24 && \\
%$\dim$(image) &   && 0 && 12 && 49 && 77 && 70 && 24&& \\
Betti num &   && 0 && 0 && 0 && 0 && 1 && 0 &&\\
\hline
\end{tabular}
\end{center} 

The table above says that $\ds  \dim \HGF{7}{2}{0}{16}=1$.  
The bases of  $\ds  \CGF{m}{2}{0}{16}$ ($m=6,7,8$) are found
either 
in Appendix 
\ref{Ap:A}, 
\ref{Ap:B} and  
\ref{Ap:C} 
or on   
\cite{mikami:URL}.

Concerning  
$\ds  \HGF{9}{2}{ }{14}$, we refer to 
\cite{metoki:shinya}, where we see the complete data. But, the notation there
is different from ours, and it seems hard to find an applicable translation rule.  
So we need to get suitable bases for our notation and 
begin searching bases without the
$k_1=0$ condition at the beginning and we get the complete bases. 
In the following discussion, we only need information about 
the bases of 
$\ds \frakC^{8}$, $\ds \frakC^{9}$ and the matrix
representation $\ds  \overline{M}$ of
$\ds\mydz: \frakC^{8} \rightarrow \frakC^{9}$, where
$\ds\frakC^{k}= \CGF{k}{2}{}{14}$.    
A similar table is obtained in the case of $\ds \frakham_{2}$  
and weight 14, 
$\rank$ is the rank of 
$ \ds  \mydz : \frak{C} ^k \longrightarrow
\frak{C}^{1+k}$.  
\begin{center}
\begin{tabular}{|c|*{8}{c}c|}
\hline
$\ds  \frakham_{2}^{ }$, wt=14 
 & & $\rightarrow$ &
$\ds  \frakC^{8} $& $ \rightarrow $ &
$\ds  \frakC^{9} $& $ \rightarrow $ & 
$\ds  \frakC^{10}$ & $ \rightarrow $ & $\mathbf{0}$\\
\hline
$\dim$  &   && 232 && 113 && 25 &&    \\
rank & & 145&& 87 && 25 & & 0 &    \\
%$\dim$(ker) &   && 145 && 88 && 25 &&     \\
%$\dim$(image) &   && 145 && 87 &&25 &&     \\
Betti num &   && 0 && 1 && 0 &&    \\
\hline
\end{tabular}
\end{center}

%A way to give a proof to our Theorem is as follows: 
3 steps of our proof of Theorem are as follows: 
\begin{enumerate}
\item To find a vector $\ds \mathbf{h} \in \ker(\mydo:
\CGF{7}{2}{0}{16} \rightarrow  
\CGF{8}{2}{0}{16})$  but 
$\ds \mathbf{h} \not\in  \mydo( \CGF{6}{2}{0}{16})$.   
\item To calculate $\ds \omega\wedge\mathbf{h}$. 
\item To check whether $\ds \omega\wedge\mathbf{h} 
\in  \mydz( \CGF{8}{2}{ }{14})$ or not. 
% ,  by counting the 
%dimension of the space generated by     
%$\ds \omega\wedge\mathbf{h}$ and $\ds\mydz( \CGF{8}{2}{ }{14})$.  
\end{enumerate}

\subsection{Gr\"obner Basis theory for cohomology groups}
To complete the proof in the direction, 
we make use of the Gr\"obner Basis theory (cf.
\cite{Cox:Little:OShea}) for 
linear homogeneous polynomials.   

Suppose we have indeterminate variables 
$(y_j)$ and fix a monomial order, say $y_1 \succ  \cdots \succ  y_{\mu}$.  
If $g_i$ are   linear homogeneous polynomials of   $(y_j)$,   then we may
write $\ds [g_1,\ldots,g_{\lambda}] = [y_1,\ldots,y_{\mu}] M $ for some
matrix $M$.  
It is well-known that we can deform $M$ into the unique elementary matrix
$\ds \widehat{M}$ (sometimes called the stairs matrix in the strict sense or
row echelon matrix) by a sequence of the three kinds elementary row
operations.  It is known as the Gaussian elimination method.  The monic
Gr\"obner basis of $\ds (g_1,\ldots g_{\lambda})$, we denote as $\ds
\text{mBasis}([g_1,\ldots,g_{\lambda}],\text{Ord}_{y})$, satisfies   
$$\ds [\text{mBasis}([g_1,\ldots,g_{\lambda}],\text{Ord}_{y}),0,\ldots,0] =
[y_1,\ldots,y_{\mu}]
\widehat{M} \ . $$  
Thus, $\rank M = \rank \widehat{M}$ is equal to the cardinality of $\ds
\text{mBasis}([g_1,\ldots,g_{\lambda}],\text{Ord}_{y})$ and   
$\text{mBasis}([g_1,\ldots,g_{\lambda}],\text{Ord}_{y})$ gives a basis for the 
$\mR$-vector space generated by 
$g_1,\ldots,g_{\lambda}$. Hereafter, we use a reduced Gr\"obner basis, we
denote it by   
$\text{Basis}([g_1,\ldots,g_{\lambda}],\text{Ord}_{y})$, for which
we allow that 
each leading coefficient should not be 1. 
So, each $j$-th element of 
$\text{Basis}([g_1,\ldots,g_{\lambda}],\text{Ord}_{y})$ is a non-zero scalar
multiple of $j$-th element of 
$\text{mBasis}([g_1,\ldots,g_{\lambda}],\text{Ord}_{y})$.

%say $y_1 \succ  \cdots \succ  y_{\mu}$.  
The 
normal form of a given polynomial $h$ with respect to a Gr\"obner basis
GB 
together with a fixed monomial order, for example
$\text{NF}(h,\text{GB},\text{Ord}_{y})$, 
is the ``smallest'' remainder of $h$ modulo by the Gr\"obner basis
$\text{GB}$.  Again, if we restrict our discussion to the linear
homogeneous polynomials, then  
$\text{NF}(h,\text{GB},\text{Ord}_{y}) =0 $ is equivalent to $\ds 
h \in $ the linear space spanned by {GB}.

We recall a key technique involving the Gr\"obner Basis theory into cohomology group
theory.  
Let $X$, $Y$ and $Z$ be finite dimensional vector spaces with bases 
$\ds\{\mathbf{q}_i\}_{i=1}^{\lambda}$, 
$\ds\{\mathbf{w}_j\}_{j=1}^{\mu}$  and 
$\ds\{\mathbf{r}_k\}_{k=1}^{\nu}$ respectively.  
Assume that there are linear maps 
$\ds g : X \rightarrow Y$ and 
$\ds f : Y \rightarrow Z$ 
%($f \circ g=0$ if we deal with cochain complexes). 
whose matrix representations are $M$ and $N$ respectively: i.e., 
\begin{align}
[g(\mathbf{q}_1),g(\mathbf{q}_2), \ldots, g(\mathbf{q}_{\lambda})]
= & [\mathbf{w}_1,\mathbf{w}_2, 
\ldots, \mathbf{w}_{\mu}] M\label{mat:rep:g}\\ \noalign{and}  
[f(\mathbf{w}_1),f(\mathbf{w}_2), \ldots, f(\mathbf{w}_{\mu})] =
&[\mathbf{r}_1,\mathbf{r}_2, \ldots ,\mathbf{r}_{\nu}] N 
\end{align}

In the right-hand side of (\ref{mat:rep:g}), 
we replace $\ds \mathbf{w}_j$ by indeterminate variable $\ds y_j$
($j=1,2,\ldots,\mu$), and get a
set of linear homogeneous polynomials $\ds [y_1,y_2,\ldots,y_{\mu}] M $.  
Denote them by $\ds [g_1(y),g_2(y),\ldots, g_{\lambda}(y)] $, i.e.,   
$\ds [g_1(y),g_2(y),\ldots,g_{\lambda}(y)] =[y_1,y_2,\ldots,y_{\mu}] M $.  

\begin{kmProp}[\cite{GB:J:Merker}]\rm \label{Merker:one}
        $GB_{e} = \text{Basis}([g_1,g_2,\ldots,g_{\lambda}],
        [y_1,y_2,\ldots,y_{\mu}],\text{Ord}_{y})$ gives a basis of $g(X)$ in the sense that   
$\ds  \{ \varphi(\mathbf{w}_1,\mathbf{w}_2,\ldots,
        \mathbf{w}_{\mu}) \mid \varphi\in 
GB_{e} \}$ forms a
basis of $g(X)$ and $\ds\rank(g) = \#( GB_{e})$. 
\end{kmProp}

We study $\ds f^{-1}(0) = \ker(f:Y \rightarrow Z)$. 
Since $\ds \langle f(\mathbf{u}), \sigma \rangle = 
\langle \mathbf{u}, f^{*}(\sigma) \rangle$ for  $\mathbf{u}\in Y, \sigma \in
Z^{*}(=\text{the dual space of } Z)$, $\ds f^{-1}(0) = 
\text{Im}(f^{*}) ^{0}$, the annihilator subspace of $\ds \text{Im}(f^{*})$.      
By Proposition \ref{Merker:one},   
we know well about   
$\ds \text{Im}(f^{*})$ by the Gr\"obner Basis theory as follows:  
Since $N$ is the matrix representation of $f$,  
$\ds {}^t N$ is a matrix representation of $\ds f^{*}$.    
We put 
$\ds   
[c_1,c_2,\ldots, c_{\mu} ] ({}^{t}N) $ by $\ds [f_1,f_2,\ldots,f_{\nu}]$.  
Fix the monomial order $\ds\text{Ord}_{c}$ of $\ds (c_j)$ by  
$c_1 \succ  \cdots \succ  c_{\mu}$.  
We get the Gr\"obner basis $\ds GB_{tr(f)}=\text{Basis}( 
[f_1,f_2,\ldots, f_{\nu} ],\;  \text{Ord}_{c})$, which gives a basis of 
$\ds \text{Im}(f^{*})$.   

Consider the polynomial $\ds h = \sum_{j=1}^{\mu} c_j y_j$, 
where $\ds  \{ y_1, \ldots, y_{\mu}\}$ are the other   
auxiliary variables (which appear for the linear map $g$). 

\begin{kmProp}[\cite{GB:J:Merker}]\rm 
        The normal form $\text{NF}(h, GB_{tr(f)}, \text{Ord}_{c})$ of $h$ 
        is written as $\ds  \sum_{j=1}^{\mu} c_j
\tilde{f}_j (y)$ where $\ds   \tilde{f}_j (y)$ is linear in $\ds  \{ y_1,
\ldots, y_{\mu}\}$. 

Let 
$\ds GB_{k} = \text{Basis}( [
 \tilde{f}_1 (y), 
 \tilde{f}_2 (y), \ldots, 
 \tilde{f}_{\mu} (y)], \text{Ord}_{y})$.   
Then $GB_{k}$ gives a basis of 
the kernel space $\ds f^{-1}(0) = \ker(f)$, and  
the cardinality of $\ds GB_{k}$ is $\ds  \dim \ker(f)$.  
\end{kmProp}

Now assume that $f \circ g=0$.   
We use the Gr\"obner bases $\ds GB_{e}$ of $g$, and $\ds GB_{k}$ of $\ker(f)$
above, then we have the following.    
\begin{kmProp}[\cite{GB:J:Merker}]\rm 
        The quotient space $\ds \ker(f: Y \rightarrow Z)/\text{Im}(g:X \rightarrow Y)$ is
        equipped with the basis 
$$ GB_{k/e} = \text{Basis}( [ \text{NF}( \varphi, GB_{e}, \text{Ord}_{y}) \mid 
\varphi\in GB_{k}], \text{Ord}_{y})$$
In particular, $\ds \dim 
\left(\ker(f: Y \rightarrow Z)/\text{Im}(g:X \rightarrow Y) \right) = \#(GB_{k/e})$.   
\end{kmProp}

\begin{kmRemark}\label{rem:stress} \rm 
In the way described above consisting of three steps, 
there is some ambiguity in choosing an element $\mathbf{h}$.   
But, if we use the Gr\"obner Basis theory, we can avoid this ambiguity.  This is 
a main reason why we use the Gr\"obner Basis theory here.  
It is hard to handle big matrices, but it is easy to deal with polynomials.
This is the second small reason. 

The last reason we use the Gr\"obner Basis theory is that 
Gr\"obner Basis packages 
in such symbol calculus softwares as Maple, Mathematica, Risa/Asir (this
is freeware) and so on, become more and more reliable and faster.   

Our calculation of Gr\"obner Bases or normal forms is 
assisted by symbol calculus software Maple.  
The author also has a proof by the Gr\"obner Basis theory to
Theorem 1.1 with the 
assistance of Maple in \cite{KM:Morita:another}.   
It will be available for reading drafts of it on   
\cite{mikami:URL} 
``A proof to Kotschick-Morita Theorem for G-K-F class''.

Risa/Asir is popular among Japanese mathematicians because it is
bundled in the Math Libre Disk which is distributed at annual
meetings of the Mathematical Society of Japan.  We put the source
code and output of our computer argument for Risa/Asir on
\cite{mikami:URL} or Appendices in \cite{KM:Morita:another}. Here,
proof to the Kotschick-Morita theorem by Risa/Asir can be seen.
You can also compare the two kinds of results calculated by Maple
and Risa/Asir, and see that the final normal forms are the same,
up to non-zero scalar multiples.  

Even in the classical linear algebra argument or the Gr\"obner
Basis argument, our discussion is based on matrix representations
of the two coboundary operators.   We stress that everything
starts from the concrete bases of cochain complexes.    
\end{kmRemark}

\subsection{Selecting a generator $h$ of $\ds\HGF{7}{2}{0}{16}$}
\label{ss:type1}
As mentioned in Remark \ref{rem:stress}, the existence of 
concrete bases of our cochain complexes is important. Actually, we got them
and can handle them, but as shown in the table above, the dimensions are 
large; for example  $\ds\dim\, C^{6}=147$, $\ds\dim\, C^{7}=95$ and 
$\ds\dim\, C^{8}=24$, where $\ds C^{k}= \CGF{k}{2}{0}{16}$.    
%It is difficult to show them all in this paper.  
Here we only show
several elements, whose number of terms of summation is smaller.    
The entire data of our concrete bases is found either in Appendix
\ref{Ap:A}, \ref{Ap:B} and \ref{Ap:C} or on 
\cite{mikami:URL}.  
\newcommand{\mywdg}{ }
The smallest element of our basis of $\ds  C^{6}$ is next, and 
consists of 28 terms: 
\begin{align*} 
 & \mathbf{q}_{142}\\
=&  -\frac{8}{3} \nz{0}{4}\km \nz{1}{4}\km \nz{2}{4}\km \nz{3}{4}\km \nz{2}{5}\km \nz{6}{7}
-\nz{0}{4}\km \nz{1}{4}\km \nz{3}{4}\km \nz{4}{4}\km \nz{1}{5}\km \nz{5}{7}
+\frac{1}{6} \nz{0}{4}\km \nz{1}{4}\km \nz{2}{4}\km \nz{4}{4}\km \nz{1}{5}\km \nz{6}{7}
-\frac{1}{6} \nz{0}{4}\km \nz{1}{4}\km \nz{2}{4}\km \nz{4}{4}\km \nz{0}{5}\km \nz{7}{7}
\\&
+\frac{2}{3} \nz{1}{4}\km \nz{2}{4}\km \nz{3}{4}\km \nz{4}{4}\km \nz{4}{5}\km \nz{0}{7}
-\frac{1}{2} \nz{0}{4}\km \nz{2}{4}\km \nz{3}{4}\km \nz{4}{4}\km \nz{0}{5}\km \nz{5}{7}
-\frac{8}{3} \nz{1}{4}\km \nz{2}{4}\km \nz{3}{4}\km \nz{4}{4}\km \nz{1}{5}\km \nz{3}{7}
+\nz{0}{4}\km \nz{2}{4}\km \nz{3}{4}\km \nz{4}{4}\km \nz{3}{5}\km \nz{2}{7}
\\&
-\frac{7}{3} \nz{0}{4}\km \nz{1}{4}\km \nz{2}{4}\km \nz{4}{4}\km \nz{3}{5}\km \nz{4}{7} 
+\frac{1}{6} \nz{0}{4}\km \nz{2}{4}\km \nz{3}{4}\km \nz{4}{4}\km \nz{4}{5}\km \nz{1}{7}
+\frac{1}{3} \nz{0}{4}\km \nz{1}{4}\km \nz{3}{4}\km \nz{4}{4}\km \nz{0}{5}\km \nz{6}{7}
+\frac{2}{3} \nz{0}{4}\km \nz{1}{4}\km \nz{2}{4}\km \nz{3}{4}\km \nz{1}{5}\km \nz{7}{7} 
\\&
+\frac{2}{3} \nz{0}{4}\km \nz{1}{4}\km \nz{3}{4}\km \nz{4}{4}\km \nz{3}{5}\km \nz{3}{7}
-\frac{8}{3} \nz{0}{4}\km \nz{1}{4}\km \nz{2}{4}\km \nz{3}{4}\km \nz{4}{5}\km \nz{4}{7}
-\frac{8}{3} \nz{1}{4}\km \nz{2}{4}\km \nz{3}{4}\km \nz{4}{4}\km \nz{3}{5}\km \nz{1}{7} 
-\frac{7}{3} \nz{0}{4}\km \nz{2}{4}\km \nz{3}{4}\km \nz{4}{4}\km \nz{2}{5}\km \nz{3}{7}
\\&
+\frac{11}{6} \nz{0}{4}\km \nz{1}{4}\km \nz{2}{4}\km \nz{4}{4}\km \nz{4}{5}\km \nz{3}{7}
+\frac{2}{3} \nz{0}{4}\km \nz{1}{4}\km \nz{2}{4}\km \nz{3}{4}\km \nz{5}{5}\km \nz{3}{7} 
+\nz{0}{4}\km \nz{1}{4}\km \nz{2}{4}\km \nz{4}{4}\km \nz{2}{5}\km \nz{5}{7}
+\frac{1}{3} \nz{0}{4}\km \nz{1}{4}\km \nz{3}{4}\km \nz{4}{4}\km \nz{5}{5}\km \nz{1}{7}
\\&
+\frac{11}{6} \nz{0}{4}\km \nz{2}{4}\km \nz{3}{4}\km \nz{4}{4}\km \nz{1}{5}\km \nz{4}{7}
-\frac{1}{2} \nz{0}{4}\km \nz{1}{4}\km \nz{2}{4}\km \nz{4}{4}\km \nz{5}{5}\km \nz{2}{7}
+\frac{2}{3} \nz{1}{4}\km \nz{2}{4}\km \nz{3}{4}\km \nz{4}{4}\km \nz{0}{5}\km \nz{4}{7}
+4 \nz{1}{4}\km \nz{2}{4}\km \nz{3}{4}\km \nz{4}{4}\km \nz{2}{5}\km \nz{2}{7}
\\&
-\frac{1}{6} \nz{0}{4}\km \nz{2}{4}\km \nz{3}{4}\km \nz{4}{4}\km \nz{5}{5}\km \nz{0}{7}
+4 \nz{0}{4}\km \nz{1}{4}\km \nz{2}{4}\km \nz{3}{4}\km \nz{3}{5}\km \nz{5}{7}
-\nz{0}{4}\km \nz{1}{4}\km \nz{3}{4}\km \nz{4}{4}\km \nz{4}{5}\km \nz{2}{7}
+\frac{2}{3} \nz{0}{4}\km \nz{1}{4}\km \nz{3}{4}\km
\nz{4}{4}\km \nz{2}{5}\km \nz{4}{7} 
\end{align*} 
where we omit the symbol $\wedge$ of wedge product. 
The two small-size elements of our basis of $\ds  C^{7}$ are the
following: 
\begin{align*} 
\mathbf{w}_{6}=&\nz{0}{3}\km\nz{1}{3}\km\nz{2}{3}\km\nz{3}{3}\km\nz{0}{6}\km
\nz{3}{6}\km \nz{6}{6} -3 \nz{0}{3}\km \nz{1}{3}\km \nz{2}{3}\km
\nz{3}{3}\km \nz{0}{6}\km \nz{4}{6}\km \nz{5}{6} -3 \nz{0}{3}\km
\nz{1}{3}\km \nz{2}{3}\km \nz{3}{3}\km \nz{1}{6}\km \nz{2}{6}\km
\nz{6}{6} \\& +6 \nz{0}{3}\km \nz{1}{3}\km \nz{2}{3}\km
\nz{3}{3}\km \nz{1}{6}\km \nz{3}{6}\km \nz{5}{6} -15 \nz{0}{3}\km
\nz{1}{3}\km \nz{2}{3}\km \nz{3}{3}\km \nz{2}{6}\km \nz{3}{6}\km
\nz{4}{6} \\ \noalign{and} \mathbf{w}_{95}=&
\nz{0}{4}\km\nz{1}{4}\km\nz{2}{4}\km\nz{3}{4}\km\nz{4}{4}\km\nz{0}{5}\km\nz{5}{5}
-5\nz{0}{4}\km\nz{1}{4}\km\nz{2}{4}\km\nz{3}{4}\km\nz{4}{4}\km\nz{1}{5}\km\nz{4}{5}
+10\nz{0}{4}\km\nz{1}{4}\km\nz{2}{4}\km\nz{3}{4}\km\nz{4}{4}\km\nz{2}{5}\km\nz{3}{5}\ .  
\end{align*} 
We pick up the smallest element of our basis of $\ds  C^{8}$: 
\begin{align*} 
\mathbf{r}_{7}  =&  \nz{0}{3}\km \nz{1}{3}\km \nz{2}{3}\km \nz{3}{3}\km \nz{0}{5}\km \nz{1}{5}\km \nz{4}{5}\km \nz{5}{5}
-2 \nz{0}{3}\km \nz{1}{3}\km \nz{2}{3}\km \nz{3}{3}\km \nz{0}{5}\km \nz{2}{5}\km \nz{3}{5}\km \nz{5}{5} 
+10 \nz{0}{3}\km \nz{1}{3}\km \nz{2}{3}\km \nz{3}{3}\km \nz{1}{5}\km \nz{2}{5}\km \nz{3}{5}\km \nz{4}{5}
\ . 
\end{align*} 
We have the matrix representations $M$ of $\ds \mydo: C^6\rightarrow C^7$  and $N$ of $\ds \mydo: C^7\rightarrow C^8$  
with respect to the
above bases.  Since the size of matrix $M$ is (95, 147) and that
of $N$ is (24,95), we will not show them here.  We only need a
generator of $\ds  \HGF{7}{2}{0}{16}$ by the Gr\"obner Basis
theory where 
we write down our linear functions $\ds  \{g_i\}$ corresponding to
$\ds  \mydo: C^{6} \rightarrow C^{7}$ and linear functions
$\{f_j\}$, giving the kernel condition for $\ds  \mydo: C^{7}
\rightarrow C^{8}$ as follows.  $$ [g_1,\ldots, g_{147}] = [
y_1,\ldots,y_{95}] M, \quad [f_1,\ldots, f_{24}] = [
c_1,\ldots,c_{95}] {}^t N $$ 
where the precise complete data are 
found either 
in Appendix \ref{Ap:D} and \ref{Ap:E} or 
on \cite{mikami:URL}.  
Here, we show a few terms as examples:  
\begin{align*}
g_{1}  =&  176 y_{1}-\frac{1036}{3} y_{8}+\frac{632}{3}
y_{9}+\frac{544}{3} y_{10}-60 y_{11}-22 y_{12}+152 y_{13}- 802
y_{21}\\&
+531 y_{22}+590 y_{23} -\frac{1625}{3} y_{24}+\frac{292}{3}
y_{25}+60 y_{26}-\frac{1595}{3} y_{27}+ 805 y_{28}\\&
-90 y_{29}+48
y_{30}-108 y_{31}-144 y_{32}-306 y_{33} +144 y_{34}+450 y_{35}\\&
+ 36 y_{36}+168 y_{37}\\ 
\vdots &\\
g_{147}  =&  \frac{5}{2}
y_{60}-\frac{7}{2} y_{61}+\frac{11}{10} y_{62}+\frac{11}{6}
y_{63}-\frac{21}{2} y_{65}+\frac{33}{10} y_{66}+\frac{15}{2}
y_{67}-\frac{1}{10} y_{68}\\&
+\frac{11}{2} y_{69}-3 y_{79}-\frac{1}{2}
y_{80} +\frac{3}{2} y_{86}+\frac{95}{12} y_{87}-\frac{17}{6} y_{88 }+2
y_{89}-\frac{209}{30} y_{90}\\&
+\frac{23}{10} y_{91}+\frac{6}{25}
y_{92}+\frac{5}{2} y_{93}-\frac{35}{2} y_{94}-\frac{133}{30} y_{95 }\\ 
\end{align*}
and 
\begin{align*}
f_{1}=& 
-\frac{55}{4} c_{1}+25 c_{3}+8 c_{5}-\frac{475}{54}
c_{8}+\frac{145}{9} c_{9}-\frac{995}{54} c_{10}+\frac{70}{3}
c_{11}+\frac{1700}{81} c_{12}\\&
-\frac{10}{81} c_{13}-\frac{41}{9}
c_{14}-\frac{215}{18} c_{15}
-\frac{425}{18} c_{16}+\frac{425}{36}
c_{17}+\frac{35}{9} c_{18}-\frac{59}{9} c_{19}+\frac{92}{9}
c_{20}\\&
+\frac{75}{32} c_{21}+\frac{85}{48} c_{22}+\frac{33}{16}
c_{23}+\frac{139}{64} c_{24}+\frac{65}{64} c_{25}
+\frac{75}{32}
c_{26}
-\frac{221}{64} c_{27}+\frac{1}{24} c_{28}\\&
-\frac{65}{12}
c_{42}+\frac{35}{4} c_{43}+\frac{95}{6} c_{44}+\frac{53}{4}
c_{45}+\frac{10}{3} c_{46}+\frac{13}{4} c_{47}+2 c_{48}-\frac{9}{2}
c_{49}-\frac{3}{2} c_{50}\\ 
\vdots &\\
f_{24}=& -15 c_{41}-10 c_{42}+30 c_{43}+35 c_{44}+3 c_{45}+40 c_{46}+18
c_{47}-3 c_{48}\\&
-\frac{301839}{740} c_{59}+\frac{256839}{740}
c_{60}
+\frac{1094769}{740} c_{61}+\frac{73128}{37}
c_{62}+\frac{174258}{185} c_{63}\\&
+\frac{848394}{185} c_{64}-\frac{435089}{740} c_{65}-\frac{105235}{111}
c_{66}-\frac{9623}{148} c_{67}
-\frac{28657}{37} c_{68} \\&
-\frac{70569}{185}
c_{69}-\frac{150326}{111} c_{70}+\frac{35965}{111}
c_{71}+\frac{4484}{185} c_{72}-\frac{3827}{37} c_{73}-\frac{52105}{74}
c_{74}\\&
+\frac{68225}{148} c_{75}-\frac{2556}{185}
c_{76}+\frac{31601}{370} c_{77}+\frac{37535}{148}
c_{78}-\frac{7439}{185} c_{79}+\frac{15139}{37} c_{80}\\&
+\frac{10657}{148}
c_{81}+\frac{56}{3} c_{86}
+\frac{53}{3} c_{87}+\frac{193}{3}
c_{88}+\frac{76}{3} c_{89}+7 c_{90}+4 c_{91}+\frac{85}{6}
c_{92}\\&
-\frac{71}{6} c_{93}-3 c_{94}+15 c_{95} \ .  
\end{align*}
The Gr\"obner basis $GB_{e}$ of $\ds  \{g_1,\ldots, g_{147}\}$
consists of 70 elements as expected.  
The whole Gr\"obner basis $GB_{e}$ is stored in Appendix \ref{Ap:I}. 
%1 2 3 4 5 6 7 8 9 0 1 2 3 4 5 6
% We only show the first and last element.  
The first element of sorted $GB_{e}$ is 
\begin{align*} 
%(1)
&  
446227638468 y_{75}-258371100400 y_{76}+2677414594200 y_{77}
- 2808720072600 y_{78}\\& 
+483892450500 y_{79}+838357655220 y_{80} -871685530860 y_{81}+
1892343009627 y_{82} \\& 
-2525687071848 y_{83} -861370434243 y_{84}-625187443152 y_{85}
- 6093198421500 y_{86}\\& 
-4546246681400 y_{87} +2813196475270 y_{88}-2152132471560 y_{89}
+15133158761840 y_{90}\\& 
-9561265966665 y_{91} -2198954966322 y_{92}+9680559087150 y_{93}
+3770983597200 y_{94}\\&
+11367701561860 y_{95} 
\, , 
\end{align*}
%& \vdots\\
and the last element of sorted $GB_{e}$ is 
%HHH
\begin{align*} &  
7228887743181600 y_{1} + 26505921724999200 y_{47}-8835307241666400
y_{50} \\& 
-40863295992707100 y_{51} - 23594575360435200 y_{76} +141145959892004100
y_{77}\\&
-152234378969531760 y_{78} + 12641923750905900 y_{79}
+103786265245653540 y_{80}\\& 
-230406289763969880 y_{81}+ 55341457461003915 y_{82} -182139922299308040 y_{83}\\&
-331644059730112995 y_{84} - 37012865309023290
y_{85}-467330302598009400 y_{86}\\& 
-327107341696261500 y_{87} + 182002246883284410 y_{88}
-27682638636383280 y_{89}\\&
+811513254542111160 y_{90} - 530692255768745745 y_{91}
-216500557914020694 y_{92}\\& 
+752677963524690150 y_{93}- 117774780478277550 y_{94} +796136446567690060 y_{95} 
\, .  
\end{align*} 
The Gr\"obner basis $GB_{k}$ corresponding to the kernel space 
defined by 
$\ds  \{f_1,\ldots, f_{24}\}$ consists of 71 elements. 
The whole 
Gr\"obner basis $GB_{k}$ is stored in Appendix \ref{Ap:J} 
or on \cite{mikami:URL}.  
%We show the first and last element below.  
The first element of sorted $GB_{k}$ is: 
\begin{align*} 
&  2027141067600 y_{76}+6871115344500 y_{77}-8293793595120 y_{78} +1593871052400 y_{79}\\& 
+3342315930030 y_{80}+2188718191440 y_{81} +6047944018587 y_{82}\\& 
-7911486513648 y_{83} +1366183084077 y_{84} -1206881491512 y_{85}\\& 
-10895090886900 y_{86}-9572836551300 y_{87} +1269138903120 y_{88}\\& 
+3867959161440 y_{89}+28054435525860 y_{90} -23511502274085 y_{91}\\& 
-7468180349703 y_{92}+2799062316375 y_{93} +17517045194250 y_{94}\\& 
+24368226519980 y_{95}\, , 
\end{align*} 
%&\vdots\\ 
and the last element of sorted $GB_{k}$ is: 
\begin{align*}
&  368571103200 y_{1}+1351427378400 y_{47}-450475792800 y_{50} -2083450541700 y_{51}\\& 
+11274054788700 y_{77}-12683683914000 y_{78} +1590428838900 y_{79}\\& 
+7275101984700 y_{80} -10448587809000 y_{81} +6410733790653 y_{82}\\& 
-13981567014072 y_{83}-16098409761957 y_{84} -2603346695478 y_{85}\\& 
-30292840571400 y_{86}-22358770705700 y_{87} +10032702042550 y_{88}\\& 
+883984609200 y_{89}+58024375642760 y_{90} -41010508144455 y_{91}\\& 
-15470397459450 y_{92} +40037017987050 y_{93} +4390494333150 y_{94}\\& 
+55052825955540 y_{95}\, . 
\end{align*}

The Gr\"obner basis corresponding to 
$\ds  \HGF{7}{2}{0}{16}$ is 
\begin{align} 
h = & 2027141067600 y_{76}+6871115344500 y_{77}-8293793595120 y_{78} \notag\\& 
+ 1593871052400 y_{79} +3342315930030 y_{80}+2188718191440 y_{81} \notag\\& 
+6047944018587 y_{82} -7911486513648 y_{83} +1366183084077 y_{84} \notag\\& 
-1206881491512 y_{85}-10895090886900 y_{86} -9572836551300 y_{87} \label{cohom:gen}\\& 
+1269138903120 y_{88}+3867959161440 y_{89}+ 28054435525860 y_{90} \notag\\& 
-23511502274085 y_{91} -7468180349703 y_{92}+2799062316375 y_{93}\notag\\& 
+17517045194250 y_{94} +24368226519980 y_{95}\ .  \notag 
\end{align}
\begin{kmRemark}\rm 
So far, all the results above, using Maple software, are also
calculated by Risa/Asir;  the outputs are found either in Appendix
\ref{Ap:L} and \ref{Ap:M}  or on  \cite{mikami:URL} by the
following title:  Results by Risa/Asir for wt=16 type1
\verb+C^{6}+\verb+->C^{7}+\verb+->C^{8}+. 

We have the generator of $\ds  \HGF{7}{2}{0}{16}$ by two methods.   
One is $h$
above by Maple.  We can see that the generator derived by Risa/Asir is $-h$; 
namely, the negative sign is the only difference.   
\end{kmRemark}

\subsection{Gr\"obner basis of $\ds\mydz ( \CGF{8}{2}{}{14})$}
\label{ss:type0} 
Next, we only need information about 
$\ds  \mydz: \frakC^{8} \rightarrow \frakC^{9}$; namely, the
bases of 
$\ds \frakC^{8}$, $\ds \frakC^{9}$ and the matrix
representation $\ds  \overline{M}$ of
$\ds  \mydz: \frakC^{8} \rightarrow \frakC^{9}$.  
These are found either in Appendix \ref{Ap:F}, \ref{Ap:G} and
\ref{Ap:H} or on \cite{mikami:URL}.  
Below we only show one of them: 
One of the 232 elements of our basis of $\ds  \frakC^{8}$
is:
\begin{align*}
\overline{\mathbf{q}}_{231}
=&  -\frac{1}{2} \nz{0}{3}\km\nz{1}{3}\km \nz{2}{3}\km \nz{3}{3}\km \nz{0}{4}\km \nz{1}{4}\km \nz{2}{4}\km \nz{6}{6}
+ \nz{0}{3}\km \nz{1}{3}\km \nz{2}{3}\km \nz{3}{3}\km \nz{0}{4}\km \nz{1}{4}\km \nz{3}{4}\km \nz{5}{6}
-\frac{1}{2} \nz{0}{3}\km \nz{1}{3}\km \nz{2}{3}\km \nz{3}{3}\km \nz{0}{4}\km \nz{1}{4}\km \nz{4}{4}\km \nz{4}{6}
\\&
-\frac{3}{2} \nz{0}{3}\km \nz{1}{3}\km \nz{2}{3}\km \nz{3}{3}\km \nz{0}{4}\km \nz{2}{4}\km \nz{3}{4}\km \nz{4}{6}
+\nz{0}{3}\km \nz{1}{3}\km \nz{2}{3}\km \nz{3}{3}\km \nz{0}{4}\km \nz{2}{4}\km \nz{4}{4}\km \nz{3}{6}
-\frac{1}{2} \nz{0}{3}\km \nz{1}{3}\km \nz{2}{3}\km \nz{3}{3}\km \nz{0}{4}\km \nz{3}{4}\km \nz{4}{4}\km \nz{2}{6}
\\&
+2 \nz{0}{3}\km \nz{1}{3}\km \nz{2}{3}\km \nz{3}{3}\km \nz{1}{4}\km \nz{2}{4}\km \nz{3}{4}\km \nz{3}{6}
-\frac{3}{2}\nz{0}{3}\km \nz{1}{3}\km \nz{2}{3}\km \nz{3}{3}\km \nz{1}{4}\km \nz{2}{4}\km \nz{4}{4}\km \nz{2}{6}
+\nz{0}{3}\km \nz{1}{3}\km \nz{2}{3}\km \nz{3}{3}\km \nz{1}{4}\km \nz{3}{4}\km \nz{4}{4}\km \nz{1}{6}
\\&
-\frac{1}{2} \nz{0}{3}\km \nz{1}{3}\km \nz{2}{3}\km \nz{3}{3}\km \nz{2}{4}\km \nz{3}{4}\km \nz{4}{4}\km \nz{0}{6} 
\end{align*} 
and, one of the 113 elements of our basis of $\ds  \frakC^{9}$ is:
\begin{align*}
\overline{\mathbf{w}}_{95} 
=&-\frac{1}{5}\nz{0}{1}\km\nz{1}{1}\km\nz{0}{4}\km\nz{1}{4}\km\nz{2}{4}\km\nz{3}{4}\km\nz{4}{4}\km\nz{0}{5}\km\nz{5}{5}
+\nz{0}{1}\km\nz{1}{1}\km\nz{0}{4}\km\nz{1}{4}\km\nz{2}{4}\km\nz{3}{4}\km\nz{4}{4}\km\nz{1}{5}\km\nz{4}{5}
-2 \nz{0}{1}\km\nz{1}{1}\km\nz{0}{4}\km\nz{1}{4}\km\nz{2}{4}\km\nz{3}{4}\km\nz{4}{4}\km\nz{2}{5}\km\nz{3}{5}
\, . 
\end{align*} 

The matrix $\overline{M}$ of 
$\ds  \mydz: \frakC^{8} \rightarrow \frakC^{9}$ is of size
(113,232) and 
the linear functions $\ds  \{\overline{g}_i\}$ corresponding to 
$\ds  \mydz: \frakC^{8} \rightarrow \frakC^{9}$ are given by 
$$ [\overline{g}_1,\ldots \overline{g}_{232}] = [ y_1,\ldots,y_{113}]
\overline{M}\, .$$ 

% \overline{g}_{1} =& 
% \vdots
% \overline{g}_{232} =& 

We will continue the same discussion as in the subsection \S
\ref{ss:type1}.  
We see that 
$\ds  \rank \overline{M} = 87$ and  
the Gr\"obner basis 
$\ds  \overline{GB}_{e} $ of 
$\ds  \{\overline{g}_i \}$, which  
corresponds to $\ds  \mydz( \frakC^{8})$, 
consists of 87 elements as expected. 
The complete data of $\ds  \{\overline{g}_i \}$, in other
words, that of $\ds  \overline{M}$, and the detail of 
$\ds  \overline{GB}_{e} $  
are found either in Appendix \ref{Ap:H} and \ref{Ap:K} 
or on \cite{mikami:URL}.

\subsection{$\omega\wedge h$ is not in $\mydz(
\CGF{8}{2}{}{14})$}

We have the linear function $h$ of $(y_1,\ldots,
y_{95})$ 
in (\ref{cohom:gen}); we know that 
the cochain $\ds  \mathbf{h}(\mathbf{w}) $ is a
non-exact kernel element in $\ds  \CGF{7}{2}{0}{16}$. 
We analyze the next element $$\ds  
\omega  \wedge h(\mathbf{w}) = 
\nz{0}{1} \wedge \nz{1}{1} \wedge h(\mathbf{w}) 
 $$ by the basis of $\ds  \frakC^{9}$,
and 
we have a linear function $\ds  \overline{h}$ of
$\ds  (y_{1},\ldots,y_{113})$ satisfying 
$$\ds  
\overline{h}(\overline{\mathbf{w}}) 
= 
\omega  \wedge h(\mathbf{w}) = 
\nz{0}{1} \wedge \nz{1}{1} \wedge h(\mathbf{w}) 
$$ which is given by the following:  
\begin{align*}
\overline{h} =&
-6996191251500 y_{74}-1557312364575 y_{76}+2027141067600 y_{77}\\&
+ 6871115344500 y_{78} -8293793595120 y_{79}+1593871052400 y_{80}\\&
+3342315930030 y_{81 }+3576568317699 y_{82} -1206881491512 y_{83}\\&
-3952406350359 y_{84}-21353158325775 y_{85}-21096249215580 y_{86}\\&
-9572836551300 y_{87}+3867959161440 y_{88}+ 10699190322480 y_{89}\\&
-23511502274085 y_{90} +2799062316375 y_{91}+17460883387175 y_{92}\\&
+17517045194250 y_{93}-43245161055925 y_{94} -121841132599900 y_{95}\, .
\end{align*} 
The normal form of 
$\ds  \overline{h}$ with respect to $\ds 
\overline{GB}_e$ is 
\begin{align*} 
- \frac{1}{1191}(\, &
7443523237284708 y_{82}+10932577142466 y_{83}- 2773751000717088 y_{84}\\& 
- 8746061713972800 y_{85} -93098703351771180 y_{90}\\&
+40450535427124200 y_{91} - 30933987324063320 y_{92} \\&
+24871612999110150 y_{93} +54636855766752700 y_{94}\\&
+ 201445748822724700 y_{95}+1180249792365936600 y_{112}\\&
+3540749377097809800 y_{113} 
\, )
\end{align*}
and is not zero, namely 
$\ds  
\overline{h}(\overline{\mathbf{w}}) 
= 
\omega  \wedge h(\mathbf{w}) = 
\nz{0}{1} \wedge \nz{1}{1} \wedge h(\mathbf{w}) 
$ is not exact, and our proof is complete.    
\kmqed

\begin{kmRemark}\rm 
Throughout this paper, the Gr\"obner basis and the normal form are computed by
Maple. On the other hand, the results by Risa/Asir are found
either 
in Appendix \ref{Ap:L} and \ref{Ap:M} or     
on \cite{mikami:URL}.

We denote by $\ds B_{maple}$ the 
normal form of $\ds \overline{h}$ with respect to $\ds 
\overline{GB}_e$ and  by $\ds A_{asir}$, the normal form calculated by
Risa/Asir.  The two are related as   
$$ B_{maple} = - \frac{1}{1191}\cdot
\frac{7443523237284708}{5337006161133135636} 
A_{asir}\ . $$ 
\end{kmRemark}

\paragraph{Acknowledgments}
%The author appreciates the referees' kind 
%advice and comments.  
The author expresses his sincere thanks to Professor
Tadayoshi Mizutani for his encouragement and discussion between
himself and the author.

%% \begin{figure}
%% \begin{center}
%% \includegraphics{figure.eps}
%% \caption{Caption}
%% \end{center}
%% \label{area}
%% \end{figure}

%% \begin{thebibliography}{9, 99 or Abc99}
%% \begin{thebibliography}{9}  for 1-digit labels
%% \begin{thebibliography}{99}  for 2-digit labels
%% \begin{thebibliography}{Abc}  for alphanumeric labels
\def\cprime{$'$} \def\cprime{$'$}

\vspace{1cm}

% \profile{Kentaro \textsc{Mikami}}
% {Akita University\\
% 1-1 Tegata Akita City, Japan\\
% mikami@math.akita-u.ac.jp}

%\profile{Firstname \textsc{Surname}}
%{Organization /Affiliation name \\
%Address\\
%E-mail}

In this paper, 
the dual basis of $R$-homogeneous polynomials of the Darboux
coordinates 
$x,y$ is denoted by $\ds \nz{r}{R}$ for $r=0,1,\ldots, R$. In
Maple, we use  the notation \texttt{z[r,R]} corresponding to  
$\ds\nz{r}{R}$, 
and the wedge product is given by 
\verb+'&^`(+ z[a,A], z[b,B], \ldots \verb+)+ and you will see them 
in the corresponding ancillary files 
indicated in the following appendices.

\begin{appendices}
        \renewcommand{\thesection}{\arabic{section}}
%\begin{tiny}
        \section{Basis of $\displaystyle \CGF{6}{2}{0}{16}$ } 
        \label{Ap:A}

%\renewcommand{\verbatiminput}{\kmcomment}

%    \verbatiminput{t1-16-6-basis.ttt} 
Please refer to the ancillary file \textbf{anc/t1-16-6-basis.ttt} 
whose file size is 2,998,777byte.

        \section{Basis of $\displaystyle \CGF{7}{2}{0}{16}$ } 
        \label{Ap:B}
%        \verbatiminput{t1-16-7-basis.ttt}
Please refer to the ancillary file \textbf{anc/t1-16-7-basis.ttt} 
whose file size is 1,561,456byte.

        \section{Basis of $\displaystyle \CGF{8}{2}{0}{16}$ } 
        \label{Ap:C}
%        \verbatiminput{t1-16-8-basis.ttt}
Please refer to the ancillary file \textbf{anc/t1-16-8-basis.ttt} 
whose file size is 138,153byte.  

        \section{Matrix representation of $\ds\mydo$ of  
        $\ds \CGF{6}{2}{0}{16}$ to $\ds \CGF{7}{2}{0}{16} $} 
        \label{Ap:D}
%       \verbatiminput{Mat_rep_t1-w16-6and7.ttt}
Please refer to the ancillary file \textbf{anc/Mat\_rep\_t1-w16-6and7.ttt} 
whose file size is 69,024byte.  
       
        \section{Dual representation of $\ds\mydo$ of 
        $\ds\CGF{7}{2}{0}{16}$ to $\ds\CGF{8}{2}{0}{16} $} 
        \label{Ap:E}
%       \verbatiminput{Dual_rep_t1-w16-7A8.ttt}
Please refer to the ancillary file \textbf{anc/Dual\_rep\_t1-w16-7A8.ttt} 
whose file size is 13,686byte.  

        \section{Basis of $\displaystyle \CGF{8}{2}{}{14}$ } 
        \label{Ap:F}
%        \verbatiminput{t0-14-8-basis.ttt}
Please refer to the ancillary file \textbf{anc/t0-14-8-basis.ttt} 
whose file size is 5,127,961byte.  

        \section{Basis of $\displaystyle \CGF{9}{2}{}{14}$ } 
        \label{Ap:G}
%        \verbatiminput{t0-14-9-basis.ttt}
Please refer to the ancillary file \textbf{anc/t0-14-9-basis.ttt} 
whose file size is 1,993,702byte.  

        \section{Matrix representation of $\ds\mydz$ of 
        $\ds\CGF{8}{2}{ }{14}$ to $\ds\CGF{9}{2}{ }{14} $} 
        \label{Ap:H}
%        \verbatiminput{Mat_rep_t0-w14-8A9.ttt}
Please refer to the ancillary file \textbf{anc/Mat\_rep\_t0-w14-8A9.ttt} 
whose file size is 98,530byte.  

        \section{Gr\"obner Basis of $\mydo (\CGF{6}{2}{0}{16})$} 
        \label{Ap:I}
%        \verbatiminput{GB_e_type1_16_6and7.ttt} 
Please refer to the ancillary file \textbf{anc/GB\_e\_type1\_16\_6and7.ttt}  
whose file size is 37,189byte.  

        \section{Gr\"obner Basis of\ $\ds\ker\mydo :\CGF{7}{2}{0}{16})
       \rightarrow  \CGF{8}{2}{0}{16}) $} 
        \label{Ap:J}
%        \verbatiminput{GB_k_type1_16_7and8.ttt}
Please refer to the ancillary file \textbf{anc/GB\_k\_type1\_16\_7and8.ttt} 
whose file size is 33,444byte.  

        \section{Gr\"obner Basis of $\mydz (\CGF{8}{2}{}{14})$} 
        \label{Ap:K}
%        \verbatiminput{GB_e-t0-14-9and9.ttt}
Please refer to the ancillary file \textbf{anc/t1-16-8-basis.ttt} 
whose file size is 138,153byte.  

        \section{Gr\"obner Basis of $\mydo$ for weight 16 by
        Risa/Asir} 
        \label{Ap:L}
%\end{tiny}
        In this section, we 
make use of Risa/Asir, which is another Symbol Calculus Software, and 
show the results we got by Maple and Risa/Asir 
are the same up to non-zero scalar multiples.  

We remark that we added some line breaks so that we get better look.

\paragraph{Basis of $\ds\mydo( \CGF{6}{2}{0}{16}) \subset
\CGF{7}{2}{0}{16}$:}\ \\ 
Our source file for Risa/Asir is this:
\begin{rmfamily}
\begin{verbatim}
/* On C^{6}_{wt=16} -> C^{7}_{wt=16} 
output("exact_t1_wt16_6and7.txt")$ 
*/ 
YList = base_var_list(y, 1, 95) $ 
load("Mat_16_6and7_type1.rr")$ 
GList = [G1, G2, G3, G4, G5, G6, G7, G8, G9, G10, G11, G12, G13, G14, G15,
G16, G17, G18, G19, G20, G21, G22, G23, G24, G25, G26, G27, G28, G29, G30,
G31, G32, G33, G34, G35, G36, G37, G38, G39, G40, G41, G42, G43, G44, G45,
G46, G47, G48, G49, G50, G51, G52, G53, G54, G55, G56, G57, G58, G59, G60,
G61, G62, G63, G64, G65, G66, G67, G68, G69, G70, G71, G72, G73, G74, G75,
G76, G77, G78, G79, G80, G81, G82, G83, G84, G85, G86, G87, G88, G89, G90,
G91, G92, G93, G94, G95, G96, G97, G98, G99, G100, G101, G102, G103, G104,
G105, G106, G107, G108, G109, G110, G111, G112, G113, G114, G115, G116, G117,
G118, G119, G120, G121, G122, G123, G124, G125, G126, G127, G128, G129, G130,
G131, G132, G133, G134, G135, G136, G137, G138, G139, G140, G141, G142, G143,
G144, G145, G146, G147]$ 
ord( YList ) $          GB1 = nd_gr(GList, YList ,0,0) $
GB1 = reverse(GB1);     print(["GBe",GB1])$     output()$           
end$
\end{verbatim}

\end{rmfamily}
\paragraph{A part of the output of Groebner Basis is:} 
\begin{verbatim}
GBe=[-446227638468*y75+258371100400*y76-2677414594200*y77+2808720072600*y78
-483892450500*y79-838357655220*y80+871685530860*y81-1892343009627*y82+
2525687071848*y83+861370434243*y84+625187443152*y85+6093198421500*y86+
4546246681400*y87-2813196475270*y88+2152132471560*y89-15133158761840*y90+
9561265966665*y91+2198954966322*y92-9680559087150*y93-3770983597200*y94
-11367701561860*y95,
2231138192340*y74-3236925592800*y76+25986517801200*y77-19167959235840*y78
+7589254454700*y79+12083813535660*y80-9614968130520*y81+19627910409861*y82
-29219550942744*y83-2582960358429*y84-3902066256336*y85-48010669486950*y86
-34187627245100*y87+10980485213590*y88+6293554661580*y89+123735218701880*y90
-76475315472495*y91-27626805021798*y92+56696019952650*y93+53106035062800*y94
+101136816168220*y95, 
68 more terms (sorry for skipping the rest and leaving heavy job to you) 
]$
\end{verbatim}

\paragraph{Kernel space of $\ds  \mydo : \CGF{7}{2}{0}{16}\rightarrow \CGF{8}{2}{0}{16}$:}
\ \\
Our source file for Risa/Asir is this:
\begin{verbatim}
/*
output("kerne_t1_wt16_7and8.txt")$
*/ 
WList = base_var_list( w, 1, 24) $ 
CList = base_var_list( c, 1, 95) $ 
YList = base_var_list( y, 1, 95) $ 
/* ###  On C^{7} -> C^{8} ### */ 
load("Mat_16_7and8_type1.rr")$ 
FList = [ F1, F2, F3, F4, F5, F6, F7, F8, F9, F10, F11, F12, F13, F14, F15,
F16, F17, F18, F19, F20, F21, F22, F23, F24, F25, F26, F27, F28, F29, F30,
F31, F32, F33, F34, F35, F36, F37, F38, F39, F40, F41, F42, F43, F44, F45,
F46, F47, F48, F49, F50, F51, F52, F53, F54, F55, F56, F57, F58, F59, F60,
F61, F62, F63, F64, F65, F66, F67, F68, F69, F70, F71, F72, F73, F74, F75,
F76, F77, F78, F79, F80, F81, F82, F83, F84, F85, F86, F87, F88, F89, F90,
F91, F92, F93, F94, F95]$ 
/* ################################################################# */ 
NagasaW = length(WList)$    NagasaF = length(FList)$ 
for ( Uke = [], J=1; J <= NagasaW; J++ ) { MyA = WList[J-1];  Atai = 0;
    for (K=1 ; K <= NagasaF; K++ ){ MyB = FList[K-1]; 
    Atai += diff( MyB, MyA)* CList[K-1];}
    Uke = cons(Atai, Uke ); }
print("mark A: Got dual map TF of F")$  Uke = reverse( Uke ); 

print("mark B: Got Image(TR)")$     GBadj = nd_gr( Uke, CList,0, 0); 

for (H=0, I=1; I <= NagasaF; I++){ H += CList[I-1]* YList[I-1]; } 
Hnf = p_nf(H, GBadj, CList, 0)$ 
for( MyUkez = [], T=CList; T != []; T = cdr(T)){ 
    MyA = car(T); MyV = diff( Hnf, MyA);
    MyUkez = cons( MyV, MyUkez);}

print("mark C: Annihilators of Image(TF)")$
MyUkez = reverse(MyUkez);

print("mark D: Got Ker(F)")$    GBker = nd_gr( MyUkez, YList,0, 0); 
output()$       end$
\end{verbatim} 

\paragraph{The outputs are the follows:}
\begin{verbatim}
GBk=[368571103200*y1+1351427378400*y47-450475792800*y50-2083450541700*y51
+11274054788700*y77-12683683914000*y78+1590428838900*y79+7275101984700*y80
-10448587809000*y81+6410733790653*y82-13981567014072*y83-16098409761957*y84
-2603346695478*y85-30292840571400*y86-22358770705700*y87+10032702042550*y88
+883984609200*y89+58024375642760*y90-41010508144455*y91-15470397459450*y92
+40037017987050*y93+4390494333150*y94+55052825955540*y95,

22987779706584*y2-31015258334280*y57+65679370590240*y58-2333788839033300*y77
+2847547488345120*y78-708876435791700*y79-1175639802428520*y80
-517877711109000*y81-1971062567955153*y82+2643976210072392*y83
-367831348898103*y84+382136247747228*y85+3663823436752050*y86
+2800692629390700*y87-463769141240790*y88-888541532645580*y89
-9149932266616200*y90+7727929995111435*y91+2481073046817666*y92
-1146891126002550*y93-5041594682543400*y94-8618277247810580*y95,

69 terms more
]$ 
\end{verbatim}

\paragraph{Basis of $\ds\HGF{7}{2}{0}{16}$}\ \\   
The next is a source file for Risa/Asir. GBe and GBk are data gotten above. 

\begin{verbatim}
/* 
output("Betti_t1_wt16_7and7.txt")$
*/

load("exact_t1_wt16_6and7.txt")$ 
load("kerne_t1_wt16_7and8.txt")$
YList = base_var_list(y,1, 95) $ 
for(Uke=[], T= GBk; T != []; T = cdr(T)) {
    MyA = car(T); print(MyA); Atai = p_nf( MyA, GBe, YList , 0) ;
    Uke = cons(Atai, Uke);
    }$
Uke = reverse(Uke)$     ord( YList )$
GBb = nd_gr( Uke, YList , 0,0);

\end{verbatim}

\paragraph{A basis of $\ds\HGF{7}{2}{0}{16}$ is given} 
\begin{verbatim}
GBkmode = [-2027141067600*y76-6871115344500*y77+8293793595120*y78
-1593871052400*y79-3342315930030*y80-2188718191440*y81-6047944018587*y82
+7911486513648*y83-1366183084077*y84+1206881491512*y85+10895090886900*y86
+9572836551300*y87-1269138903120*y88-3867959161440*y89-28054435525860*y90
+23511502274085*y91+7468180349703*y92-2799062316375*y93-17517045194250*y94
-24368226519980*y95]
\end{verbatim}

%\begin{tiny}
        \section{Gr\"obner Basis of $\mydz$ for weight 14 by Risa/Asir } 
        \label{Ap:M}
%\end{tiny}
        \parindent=0pt 
In this section, we make use of Risa/Asir, which is another Symbol 
Calculus Software, and show the results we got by Maple and
Risa/Asir are the same up to non-zero scalar multiples.  

We remark that we added some line breaks so that we get better look.

% \paragraph{Basis of $\ds\mydz( \CGF{8}{2}{ }{14}) \subset
% \CGF{9}{2}{ }{14}$:}\ \\ 
Our source file for Risa/Asir is this:
\begin{rmfamily}
\begin{verbatim}
/* On C^{8} -> C^{9} 
output("exact_t0_wt14_8and9.txt")$ */
YList = base_var_list(y, 1, 113)$
load("Mat_14_8and9_type0.rr")$
FList = [F1, F2, F3, F4, F5, F6, F7, F8, F9, F10, F11, F12, F13, F14, F15, F16, 
F17, F18, F19, F20, F21, F22, F23, F24, F25, F26, F27, F28, F29, F30, F31, F32, 
F33, F34, F35, F36, F37, F38, F39, F40, F41, F42, F43, F44, F45, F46, F47, F48, 
F49, F50, F51, F52, F53, F54, F55, F56, F57, F58, F59, F60, F61, F62, F63, F64, 
F65, F66, F67, F68, F69, F70, F71, F72, F73, F74, F75, F76, F77, F78, F79, F80, 
F81, F82, F83, F84, F85, F86, F87, F88, F89, F90, F91, F92, F93, F94, F95, F96, 
F97, F98, F99, F100, F101, F102, F103, F104, F105, F106, F107, F108, F109, F110, 
F111, F112, F113, F114, F115, F116, F117, F118, F119, F120, F121, F122, F123, 
F124, F125, F126, F127, F128, F129, F130, F131, F132, F133, F134, F135, F136, 
F137, F138, F139, F140, F141, F142, F143, F144, F145, F146, F147, F148, F149, 
F150, F151, F152, F153, F154, F155, F156, F157, F158, F159, F160, F161, F162, 
F163, F164, F165, F166, F167, F168, F169, F170, F171, F172, F173, F174, F175, 
F176, F177, F178, F179, F180, F181, F182, F183, F184, F185, F186, F187, F188, 
F189, F190, F191, F192, F193, F194, F195, F196, F197, F198, F199, F200, F201, 
F202, F203, F204, F205, F206, F207, F208, F209, F210, F211, F212, F213, F214, 
F215, F216, F217, F218, F219, F220, F221, F222, F223, F224, F225, F226, F227, 
F228, F229, F230, F231, F232]$  

ord( YList ) $      GB1 = nd_gr(FList, YList ,0,0) $
GB1 = reverse(GB1); print(["GBe",GB1])$
output()$           end$
\end{verbatim}

\end{rmfamily}
\paragraph{A part of the output of Groebner Basis is:} 
\begin{verbatim}
GBe = [5*y103-2*y104+y105+5*y106+y107-8*y108+4*y109-y111-14*y112-420*y113,
15*y96+10*y97-10*y99-10*y100-y104-2*y105-10*y106-2*y107+16*y108-8*y109
    +2*y111+28*y112+1540*y113,
20*y89+24*y90-24*y91+24*y92-30*y93-8*y94-80*y95-15*y105+18*y107-18*y108+24*y109
    -60*y110-224*y112,
84 more terms
]$
\end{verbatim} 

%\newpage
\paragraph{Morita conjecture for Kontsevich map 
$\ds \omega\wedge :\HGF{7}{2}{0}{16}\rightarrow \HGF{9}{2}{ }{14}$:}
\  

We have 
$h$ 
in 
$\ds \HGF{7}{2}{0}{16}$ and $\overline{h} = \omega\wedge h$ is given 
by 
\begin{verbatim}
/* This is a data of \omega \wedge h(w) in C^{9}_{wt=14) */ 
Hbar = -6996191251500*y74-1557312364575*y76+2027141067600*y77+6871115344500*y78
-8293793595120*y79+1593871052400*y80+3342315930030*y81+3576568317699*y82
-1206881491512*y83-3952406350359*y84-21353158325775*y85-21096249215580*y86
-9572836551300*y87+3867959161440*y88+10699190322480*y89-23511502274085*y90
+2799062316375*y91+17460883387175*y92+17517045194250*y93-43245161055925*y94
-121841132599900*y95$
\end{verbatim}

Our source file for Risa/Asir is this:
\begin{verbatim}
/* 
Final stage of checking hbar is in GBe or not at C^{9}_{wt=14}.  
*/

output("kekka_t0_wt14_9and9.txt")$ 
load("exact_t0_wt14_8and9.txt")$ /* GBe */ 
YList = base_var_list(y, 1, 113)$

Hbar = -6996191251500*y74-1557312364575*y76+2027141067600*y77 
+6871115344500*y78-8293793595120*y79+1593871052400*y80+3342315930030*y81
+3576568317699*y82-1206881491512*y83-3952406350359*y84-21353158325775*y85
-21096249215580*y86-9572836551300*y87+3867959161440*y88+10699190322480*y89
-23511502274085*y90+2799062316375*y91+17460883387175*y92+17517045194250*y93
-43245161055925*y94-121841132599900*y95$

ord(YList)$     Atai = p_nf(Hbar, GBe, YList, 0) ;
output()$       end$ 
\end{verbatim} 

\paragraph{The output is:}
\begin{verbatim}
5337006161133135636*y82+7838657811148122*y83-1988779467514152096*y84
-6270926248918497600*y85-66751770303219936060*y90+29003033901248051400*y91
-22179668911353400440*y92+17832946520361977550*y93+39174625584761685900*y94
+144436601905893609900*y95+846239101126376542200*y112+2538717303379129626600*y113
\end{verbatim}

\bigskip

\textbf{Remark:}
        The output just above, we denote by $\ds A_{asir}$.  
        We denote by $\ds B_{maple}$ the output in the preprint 
{\em An affirmative answer to a conjecture for Metoki class} by K.~Mikami.

You will verify that
$$ B_{maple} = - \frac{1}{1191}\cdot\frac{7443523237284708}{5337006161133135636} 
A_{asir}$$

%a := 5337006161133135636: b := 7443523237284708:

\end{appendices}

\noindent
Kentaro \textsc{Mikami}\\ 
Akita University\\
1-1 Tegata Akita City, Japan\\
mikami@math.akita-u.ac.jp
\end{spacing}
\end{document}